\magnification=1200
\baselineskip=18pt
\def\bull{\vrule height .9ex width .8ex depth -.1ex } % 
\tolerance=100000
\overfullrule=0pt

\centerline{\bf QUASI-CLASSICAL LIE-SUPER ALGEBRA}
\centerline{\bf AND LIE-SUPER TRIPLE
SYSTEMS}

\vskip .3 truein

\centerline{by}

\vskip .3 truein

\centerline{Susumu Okubo}
\centerline{Department of Physics and Astronomy}
\centerline{University of Rochester}
\centerline{Rochester, NY 14627 U.S.A.}

\vskip .3 truein

\centerline{and}

\vskip .3 truein

\centerline{Noriaki Kamiya}
\centerline{Department of Mathematics}
\centerline{Shimane University}
\centerline{Matsue, Shimane 690, Japan}

\vskip .6 truein

\centerline{\underbar{\bf Abstract}}

Notions of quasi-classical Lie-super algebra as well as Lie-super triple 
systems have been given
 and studied 
 with some examples. Its application to Yang-Baxter 
equation has also been given.

\vskip 1.3 truein

\noindent 1991 Mathematics Subject Classification Primary 17A40, 17B60

\vfill\eject

\noindent {\bf 1. \underbar{Quasi-classical Lie-super Algebra}}

\medskip

In this note, we will first introduce the notion of quasi-classical
 Lie-super algebra as well as quasi-classical Lie-super triple system with some 
examples. We will apply them to obtain some new solutions of Yang-Baxter 
equation in section 3. Algebras in this note are assumed to be finite 
dimensional over a field of characteristic not two.

Let $L$ be a Lie-super algebra, i.e. it is first a direct sum
$$L = V_0 \oplus V_1 \eqno(1.1)$$
of bosonic $(V_0)$ and fermionic $(V_1)$ spaces. We denote the grade by
$$\sigma (x) = \cases{0\ , &if $\ \ x\ \epsilon \ V_0$\cr
\noalign{\vskip 6pt}%
1\ , &if $\ \ x\ \epsilon\ V_1$ \cr}\eqno(1.2)$$
and write
$$(-1)^{\sigma (x) \sigma(y)} = (-1)^{xy}\quad . \eqno(1.3)$$
 Then, the 
Lie-product $[x,y]$ satisfies the following conditions:

\medskip

\item{(i)} $\sigma ([x,y]) = \{ \sigma (x) + \sigma (y)\} \ {\rm mod}\ 2 
$ \hfill (1.4a)
\smallskip
\item{(ii)} $[y,x] = -(-1)^{xy} [x,y]$ \hfill (1.4b)
\smallskip
\item{(iii)} $(-1)^{xy} [[x,z],y] + (-1)^{yz} [[y,x],z] + (-1)^{zx} [[z,y],
x] = 0 \quad .$ \hfill (1.4c)

\medskip

Suppose now that $L$ possesses a bilinear non-degenerate form $<.|.>$ 
satisfying conditions:

\medskip

\item{(i)} $<x|y>\ =0$, unless $\sigma (x) = \sigma (y)$
\hfill (1.5a)
\smallskip
\item{(ii)} $<y|x>\ = (-1)^{xy} <x|y>$ \hfill (1.5b)
\smallskip
\item{(iii)} $< [x,y]|z>\ =\ <x|[y,z]> \quad .$ \hfill (1.5c)

\medskip

\noindent We will then call $L$  quasi-classical. If $L$ is a simple
 Lie-super algebra with non-zero Killing form [1], we may then set
$$<x|y>\ =\ {\rm Tr} (adx\  ady)$$
where Tr  hereafter stands for the super-trace. However, the 
converse is not necessarily true as we will soon see.

Let $e_1, e_2, \dots , e_N$ with $N = \dim L$ be a basis of $L$ with
$$\eqalignno{&\sigma (e_j) = \sigma_j &(1.6a)\cr
&[e_j , e_k ] = \sum^N_{\ell = 1} C^\ell_{jk} e_\ell \quad . &(1.6b)\cr}$$
Suppose that $L$ possesses a Casimir invariant $I_2$ given by
$$\eqalignno{&I_2 = \sum^N_{j,k=1} g^{jk} e_j e_k &(1.7a)\cr
\noalign{\vskip 4pt}%
&g^{jk} = (-1)^{\sigma_j \sigma_k} g^{kj} &(1.7b)\cr
&g^{jk} = 0 \quad {\rm if}\quad \sigma_j \not= \sigma_k \quad . &(1.7c)
\cr}$$
The condition $[I_2, e_\ell] = 0$ is equivalent to the validity of 
$$\sum^N_{m=1} g^{jm} C^k_{m\ell} = \sum^N_{m=1} C^j_{\ell m} g^{mk} \quad 
. \eqno(1.8)$$

\medskip

\noindent {\bf \underbar{Proposition 1.1}}

A necessary and sufficient condition of a Lie-super algebra $L$ being
 quasi-classical is
 the existence of the Casimir invariant $I_2$ such that $g^{jk}$ 
is non-degenerate with its inverse $g_{jk}$, i.e. 
$$\eqalignno{&\sum^N_{\ell =1} g^{k \ell} g_{\ell j} = \sum^N_{\ell =1} g_{
j \ell} g^{\ell k} = \delta^k_j &(1.9a)\cr
\noalign{\vskip 4pt}%
&g_{jk} = (-1)^{\sigma_j \sigma_k} g_{kj} &(1.9b)\cr
&g_{jk} = 0 \quad {\rm unless}\quad \sigma_j = \sigma_k \quad .&(1.9c)
\cr}$$
\medskip

\noindent {\bf \underbar{Proof}}

Suppose that $L$ is quasi-classical. Setting 
$$g_{jk} =\ <e_j |e_k> \quad ,$$
it has its inverse $g^{jk}$. The relation
$$<[e_j,e_k]|e_\ell>\ =\ <e_j |[e_k,e_\ell]>$$
can  easily be shown to be equivalent to Eq. (1.8) so that $I_2$ defined by 
Eq. (1.7a) is the Casimir invariant. Conversely, let us assume that the 
Casimir invariant $I_2$ exists. We introduce the bilinear form $<.|.>$ in
 $L$ by
$$<e_j |e_k>\ = g_{jk}$$
which defines the desired bilinear non-degenerate supersymmetric form 
satisfying Eqs. (1.5).  $\quad \bull$

\medskip

\noindent \underbar{\bf Remark 1.1}

This proposition is a straightforward generalization of the result of [2].

We will now give some examples of quasi-classical Lie and Lie-super 
algebras below.

\medskip

\noindent {\bf \underbar{Example 1.1}}

Let $L = V_0 = \{ e, f, x_1, \dots , x_n , y_1 , \dots , y_n \}$ with
 $V_1 =0$. Only non-zero Lie products are assumed to be given by
$$\eqalign{[x_j ,f] &= -[f,x_j] = x_j\cr
[y_j, f] &= -[f,y_j] = -y_j\cr
[x_j, y_k] &= -[y_k ,x_j] = \delta_{jk} e\cr}$$
for $j,k = 1,2,\dots ,n$.  It is easy to verify that $L$ is a Lie algebra 
with the Casimir invariant
$$I_2 = \lambda e^2 + ef + fe - \sum^n_{j=1} (x_j y_j + y_j x_j)$$
for arbitrary constant $\lambda$. Note that $e$ 
 is a center element of $L$. We can now introduce the inner product by
$$<e|f>\ = \ <f|e>\ = 1 \quad , \quad <f|f>\ = - \lambda \quad ,$$
$$<x_j |y_k>\ =\ <y_k |x_j>\ = - \delta_{jk} \quad ,$$
while all other inner products are assumed to be zero. We can readily 
verify that $L$ is quasi-classical.

\medskip
\vfill\eject

\noindent {\bf \underbar{Example 1.2}}

Let $L = V_0 \oplus V_1$ with
$$V_0 = \{ e,f\} \quad , \quad V_1 = \{ x_1, \dots ,x_n, y_1 ,\dots ,y_n\} 
\quad ,$$
where only non-zero products are assumed to be given by
$$\eqalign{[x_j,f] &= -[f,x_j] = x_j \quad ,\cr
[y_j,f] &= -[f,y_j] = -y_j \quad ,\cr
[x_j , y_k] &= [y_k , x_j] = \epsilon_{jk} e \quad .\cr}$$
Here, $\epsilon_{jk} = - \epsilon_{kj}$ is antisymmetric with its inverse $
\epsilon^{jk}$. Especially, $n$ must now be even. The Casimir invariant is 
found to be
$$I_2 = \lambda e^2 + ef + fe + \sum^n_{j,k=1} \epsilon^{jk} \{x_j y_k -y_k 
x_j\}\quad .
$$
We introduce inner products by
$$<f|f>\ = - \lambda \quad ,\quad <e|f>\ = \ <f|e>\ = 1 \quad ,$$
$$<x_j | y_k>\ = - <y_k|x_j>\ = - \epsilon_{jk} \quad ,$$
while all other $<.|.>$ are zero. Here, $\lambda$ is again an arbitrary 
constant. We can verify that $L$ is quasi-classical.

\medskip

\noindent {\bf \underbar{Remark 1.2}}

Both examples 1.1 and 1.2 given above are not simple but solvable, since 
they satisfy the identity
$$[L,[[L,L],[L,L]]] = 0 \quad . \eqno(1.10)$$
However, they are not nilpotent since $[L,[L,L]] = [L,L] \not= 0$.
We will next give examples of nilpotent quasi-classical Lie and Lie-super
 algebras.

\medskip
\vfill\eject

\noindent {\bf \underbar{Example 1.3}}

$$L = V_0 = \{x_j, u_j, y_A, v_A, Y_{jA}\} \quad {\rm with}\quad
V_1 = 0\quad ,$$
where indices $j$ and $A$ assumes $j=1,2,\dots,n$ and $A=1,2,\dots,m$. Only 
non-zero commutators are given by
$$\eqalign{[x_j,Y_{kA}] &= -[Y_{kA},x_j] = \delta_{jk} v_A\cr
[y_A, Y_{jB}] &= -[Y_{jB}, y_A] = - \delta_{AB}u_j\cr
[x_j,y_A] &= -[y_A,x_j] = - Y_{jA}\cr}$$
for $j,k=1,2,\dots,n$ and $A,B = 1,2,\dots,m$. $L$ can be verified to be a 
Lie algebra with center elements $\{u_j,v_A\}$. The Casimir invariants is 
found to be
$$\eqalign{I_2 = &\sum^n_{j=1} (x_j u_j +u_jx_j) +
 \sum^m_{A=1} (v_A y_A + y_A v_A)\cr 
&+ \sum^n_{j=1} \sum^m_{A=1} \Lambda_{jA}\Lambda_{jA} \quad.\cr}$$
Actually, we can add bilinear terms involving center elements $u_j$ and
 $v_A$ to this expression. However, we will not do so here for simplicity. The 
corresponding inner products are calculated to be
$$\eqalign{&<Y_{jA} |Y_{kB}>\ = \delta_{jk} \delta_{AB}\cr
&<x_j |u_k>\ = \ <u_k|x_j>\ = \delta_{jk}\cr
&<v_A |y_B>\ = \ <y_B |v_A>\ = \delta_{AB}\cr}$$
while all other $<.|.>$ are zero.

\medskip

\noindent {\bf \underbar{Example 1.4}}

$$L_0 = V_0 \oplus V_1 \quad {\rm with}\quad V_0 = \{ x_j, u_j\} \quad , 
\quad V_1 =\{y_A , v_A, Y_{jA}\}$$
as in Example 1.3, except for the fact that we replace relations for $[y_A, 
Y_{jB}]$,  $<v_A |y_B>$  etc. by
$$[y_A,Y_{jB}] = [Y_{jB},y_A] = - \epsilon_{AB} u_j \quad ,
\quad <Y_{jA}|Y_{kB}>\ =\ \delta_{jk} \epsilon_{AB}\quad ,$$
$$<v_A |y_B>\ = - <y_B|v_A>\  = \epsilon_{AB}$$
for a symplectic form $\epsilon_{AB} = - \epsilon_{BA}$ with its inverse
 $\epsilon^{AB}$. The Casimir invariant $I_2$ will now be given  by
$$\eqalign{I_2 = &\sum^n_{j=1} (x_j u_j +u_j x_j) + \sum^m_{A,B=1}
 \epsilon^{AB} (v_A y_B -y_B v_A)\cr
&+ \sum^n_{j=1} \sum^m_{A,B=1} \epsilon^{AB} \Lambda_{jA}\Lambda_{jB} \quad 
. \cr}$$
\medskip

\noindent {\bf \underbar{Remark 1.3}}

Let us define $L_n (n=1,2,\dots)$ by $L_1 =L$, and $L_{n+1} = [L,L_n]$ 
inductively. If we have $L_{n+1} =0$ but $L_n \not= 0$, then we say that 
the Lie super-algebra $L$ is nilpotent of the length $n$. The examples 1.3 
and 1.4 satisfy $L_3 \not= 0$ but $L_4 =0$ so that both are nilpotent with 
length 3.

\medskip

\noindent {\bf \underbar{Remark 1.4}}

The non-degenerate bilinear form $<x|y>$ is not unique. Note that the 
examples 1.1 and 1.2 contain an arbitrary parameter $\lambda$. 
This is due to the existence of the center element $e$, as the following 
proposition will show. Some other examples of quasi-classical Lie algebras 
which are not super algebra are also found in ref. [3]. 

\medskip

\noindent {\bf \underbar{Proposition 1.2}}

Let a  Lie-super algebra $L$ possess two bilinear forms 
$<x|y>_j$ $(j=1,2)$ satisfying conditions Eqs. (1.5). Suppose that the 
adjoint representation of $L$ is irreducible i.e., that if $A               
\ \epsilon\ End\ L$
 is grade-preserving and
 satisfies $[adx, A] =0$ for all $x\ \epsilon\ L$, then 
$A = \lambda Id$ for a constant $\lambda$. Here $Id$ is the identity 
mapping in $L$.  Then, if $<x|y>_1$ is non-degenerate, we have
$$<x|y>_2\  = \lambda <x|y>_1$$
for a constant $\lambda$. We note that we need not assume the
 non-degeneracy of $<x|y>_2$.

\medskip

\noindent {\bf \underbar{Proof}}

Since $L$ is finite dimensional and since $<x|y>_1$ is assumed to be
 non-degenerate, the standard reasoning implies the existence of
 $A \ \epsilon\ End\ L$ such that
$$<x|y>_2\  =\ <Ax|y>_1 \quad .$$
Moreover, $A$ is grade-preserving, i.e. $\sigma (Ax) = \sigma(x)$. Now, the 
condition $<[y,x]|z>_j\  =$ $<y|[x,z]>_j$ is then rewritten as
$$<[A,adx]y|z>_1\  = 0$$
which leads to $[A, adx] =0$ because of the non-degeneracy of $<y|z>_1$. 
The irreducibility assumption leads to the desired result $A = \lambda Id$ 
and hence $<x|y>_2\  = \lambda <x|y>_1$.

\medskip

\noindent {\bf \underbar{Remark 1.5}}

The adjoint representation is irreducible, if $L$ is simple and,
if the underlying field is algebraically closed. $\quad \bull$

Applying a theorem due to Dieudonn\'e (see [4] p. 24) on an algebra possessing an 
associative bilinear form, we have also the following proposition.

\medskip

\noindent {\bf \underbar{Proposition 1.3}}

Suppose that we have $[B,B] \not= 0$ for every ideal $B$ of a quasi-
classical Lie-super algebra $L$. Then, $L$ is uniquely expressible as a 
direct sum
$$L = B_1 \oplus B_2 \oplus \dots \oplus B_t$$
of simple ideals $B_j$.

\medskip

\noindent {\bf 2. \underbar{Quasi-classical Lie-super Triple System}}

A $Z_2$-graded vector space $V$ is called a $\delta$ Lie-super triple 
system for $\delta = \pm1$, if it possesses a triple linear product 
$V \otimes V \otimes V \rightarrow V$ satisfying

\medskip

\item{(0)} $\sigma ([x,y,z]) = (\sigma(x) + \sigma(y) + \sigma(z))\ \ ({\rm 
mod}\ 2)$\hfill (2.1a)
\smallskip
\item{(1)} $[y,x,z] = - \delta (-1)^{xy} [x,y,z]$\hfill (2.1b)
\smallskip
\item{(2)} $(-1)^{xz} [x,y,z] + (-1)^{yx} [y,z,x] + (-1)^{zy} [z,x,y] = 0$ 
\hfill (2.1c)
\smallskip
\item{(3)} $[u,v,[x,y,z]] = [[u,v,x],y,z] + (-1)^{(u+v)x} [x,
[u,v,y],z]$

\item{ }  $\qquad\qquad\qquad\qquad +\  (-1)^{(u +v)(x+y)} [x,y,[u,v,z]] 
\quad .$ \hfill (2.1d)

\medskip

\noindent Especially, the case of $\delta =1$ defines a Lie-super triple 
system while the other case of $\delta = -1$ may be termed an\
 anti-Lie-super triple system as in [5].

Moreover, suppose that there exists a non-degenerate bilinear form $<.|.>$ 
in $V$ obeying conditions:

\medskip

\item{(1)} $<x|y>\ =0 \quad {\rm unless}\quad \sigma (x) = 
\sigma (y) $ \hfill (2.2a)
\smallskip
\item{(2)} $<y|x>\ = \delta (-1)^{xy} <x|y>$ \hfill (2.2b)
\smallskip
\item{(3)} $<[x,y,u]|v>\ = - (-1)^{(x+y)u} <u|[x,y,v]> \quad .$ \hfill 
(2.2c)

\medskip

\noindent We then call  the $\delta$ Lie-super triple system $V$  
quasi-classical.

We will first prove the following:

\medskip

\noindent {\bf \underbar{Proposition 2.1}}

Let $V$ be a  $\delta$ Lie-super triple system with a 
possible exception of the validity of Eq. (2.1d).
 Moreover assume the validity of Eq. (2.2b).
 The following 3 
conditions are then equivalent to each other:

\medskip

\item{(1)} $<[x,y,u]|v>\ = - (-1)^{(x+y)u} <u|[x,y,v]>$ \hfill
(2.3a)
\smallskip
\item{(2)} $<[x,y,u]|v>\ = -(-1)^{(u+v)y}<x|[u,v,y]>$\hfill (2.3b)
\smallskip
\item{(3)} $<x|[y,u,v]>\ = (-1)^{xy+uv} <y|[x,v,u]>
 \quad  .$ \hfill (2.3c)

\medskip

\noindent {\bf \underbar{Proof}}

\smallskip

\noindent (i) \underbar{(2) $\rightarrow$ (1)}

Letting $u \leftrightarrow v$ in (2), it gives
$$\eqalign{<[x,y,u]|v>\ &= - \delta (-1)^{uv} <[x,y,v]|u>\cr
&= -(-1)^{uv} (-1)^{(x+y+v)u} <u|[x,y,v]>\cr
&= -(-1)^{(x+y)u} <u|[x,y,v]>\cr}$$
which is (1).

\smallskip

\noindent (ii) \underbar{(3) $\rightarrow$ (2)}

$$\eqalign{<[x,y,u]|v>\ &=  \delta (-1)^{v(x+y+u)} <v|[x,y,u]>\cr
&= \delta (-1)^{v(x+y+u)} (-1)^{vx+yu} <x|[v,u,y]>\cr
&= -(-1)^{y(u+v)} <x|[u,v,y]>\cr}$$
which is (2).

\smallskip

\noindent (iii) \underbar{(2) $\rightarrow$ (3)}

Because of (i), we may assume the validity of both (1) and (2). Then
$$<u|[x,y,v]>\ = -(-1)^{(x+y)u} <[x,y,u]|v>$$
by (1). However, $<[x,y,u]|v>\ =-(-1)^{y(u+v)} <x|[u,v,y]>$ by (2). 
Combining both, we obtain
$$<u|[x,y,v]> \ = (-1)^{xu+yv} <x|[u,v,y]>\quad .$$
Interchanging $x \rightarrow y \rightarrow u \rightarrow x$, this leads to 
(3).

\smallskip

\noindent (iv) \underbar{(1) $\rightarrow$ (2)}

We first note that (1) implies
$$<[x,y,u]|v>\ = - \delta (-1)^{uv} <[x,y,v]|u> \quad . \eqno(2.4)$$
Using Eq. (2.1c), we calculate
$$\eqalign{<[x,y,u]|v>\ &= -(-1)^{(x+y)u} <u|[x,y,v]>\cr
&= (-1)^{(x+y)u} \big\{ (-1)^{xv+yx} <u|[y,v,x]>\cr
&\qquad\quad +(-1)^{xv+vy} <u|[v,x,y]>\big\}\cr
&= -(-1)^{x(u+v+y)+uv}<[y,v,u]|x>\cr
&\qquad\quad -(-1)^{v(x+y+u)+yu} <[v,x, u]|y>\quad .\cr}$$
Now, we let $u \leftrightarrow v$ and note Eq. (2.4). We calculate then
$$\eqalign{2<[x,y,u]|v>\ &=
\ <[x,y,u]|v>\ - \ \delta  (-1)^{uv} <[x,y,v]|u>\cr
&= (-1)^{x(u+v+y)+yv} \delta <(-1)^{uv} [v,y,u] +
(-1)^{yv} [y,u,v]|x> \cr
&\qquad\quad - (-1)^{v(x+y)+yu} <(-1)^{uv} [v,x,u] +
(-1)^{xv} [x,u,v]|y>\cr
&= - \delta (-1)^{(u+v)(x+y)+xy}<[u,v,y]|x>\cr
&\qquad\quad +(-1)^{(u+v)(x+y)}
 <[u,v,x]|y>\quad .\cr}$$
Now, interchanging $x \leftrightarrow u$ and $y \leftrightarrow v$ in Eq.
 (2.4), we have $<[u,v,x]|y>\ = - \delta (-1)^{xy}$
 $ <[u,v,y]|x>$ so that 
$$\eqalign{<[x,y,u]|v>\ &= -\delta (-1)^{(u+v)(x+y)+xy} <[u,v,y]|x>\cr
&= -(-1)^{(u+v)y} <x |[ u,v,y]>\cr}$$
which is (2). This completes the proof. $\quad \bull$

Next, we will define left and right multiplication operators $V \otimes V
 \rightarrow End\ V$ by
$$\eqalignno{L(x,y)z &= [x,y,z] &(2.5a)\cr
R(x,y)z &= (-1)^{z(x+y)} [z,x,y] \quad , &(2.5b)\cr}$$
and set
$$[L(u,v),R(x,y)] = L(u,v)R(x,y) - (-1)^{(x+y)(u+v)} R(x,y) L(u,v)
\eqno(2.6)$$
and similarly for $[L(u,v),\ L(x,y)]$.

\medskip

\noindent \underbar{\bf Lemma 2.1}

$$\eqalignno{&L(y,x) = -\delta (-1)^{xy} L(x,y) &(2.7a)\cr
&[L(u,v), L(x,y)] = L([u,v,x],y) + (-1)^{(u+v)x}L(x,[u,v,y])&(2.7b)\cr
&[L(u,v),R(x,y)] = R([u,v,x],y) + (-1)^{(u+v)x} R(x,[u,v,y])
\quad .&(2.7c)\cr}$$

\medskip

\noindent \underbar{\bf Proof}

Eqs. (2.7a) and (2.7b) are immediate consequences of Eqs. (2.1b) and
 (2.1d). To show Eq. (2.7c), we calculate
$$\eqalign{[L(u,v),R(x,y)]z &= (-1)^{(x+y)z} \{ [u,v,[z,x,y]]-[[u,v,z],x,y]\}\cr
&= (-1)^{(x+y +u+v)z} \{[z,[u,v,x],y] + (-1)^{x(u+v)} [z,x,[u,v,y]]\}\cr
&= R([u,v,x],y)z + (-1)^{x(u+v)} R(x,[u,v,y])z\cr}$$
which proves (2.7c). $\quad \bull$

\medskip

\noindent \underbar{\bf Proposition 2.2}

Let $V$ be a $\delta$ Lie-super triple system. If $<x|y>_1$ defined by
$$<x|y>_1 = {1 \over 2} \ {\rm Tr} \{ R(x,y) + \delta (-1)^{xy} R(y,x)\}$$
is non-degenerate, then $V$ is quasi-classical. Here Tr stands for the 
supertrace as before.

\medskip

\noindent \underbar{\bf Proof}

The conditions Eqs. (2.2a) and (2.2b) follow readily from the definition. 
Taking the supertrace of both sides, Eq. (2.7c) gives
$${\rm Tr}\ R([u,v,x],y) + (-1)^{(u+v)x}\  {\rm Tr} \ R(x,[u,v,y]) = 0$$
which leads to the validity of
$$< [u,v,x]|y>_1\  = -(-1)^{(u+v)x} <x|[u,v,y]>_1 \quad .$$

\medskip

\noindent \underbar{\bf Remark 2.1}

We can prove contrarily Tr $L(x,y) = 0$ identically.

We shall now give some examples of quasi-classical $\delta$ Lie-super 
triple system.

\medskip

\noindent \underbar{\bf Example 2.1}

Let $V$ be a $Z_2$-graded vector space with a non-degenerate bilinear form 
$<x|y>$ satisfying

\medskip

\item{(i)} $<x|y> \ =0 \quad {\rm unless}\quad \sigma (x) = \sigma (y)$
\smallskip

\item{(ii)} $<y|x>\ = \delta (-1)^{xy} < x|y> \quad .$

\medskip

\noindent Then, the triple product
$$[x,y,z] =\ <y|z>x - \delta (-1)^{xy} <x|z>y$$
defines a quasi-classical $\delta$ Lie-super triple system.

\medskip

\noindent \underbar{\bf Example 2.2}

Let $V$ be as above, and let $P\ \epsilon\ End\ V$ satisfy conditions

\medskip

\item{(i)} $\sigma (Px) = \sigma (x)$
\smallskip

\item{(ii)} $<x|Py>\ =\ <Px|y>$
\smallskip

\item{(iii)} $P^2 = c Id$

\medskip

\noindent for a constant $c$, where $Id$ stands for the identity mapping.
 The triple product defined by
$$[x,y,z] = \ <y|z>Px\ + \ <y|Pz>x - \delta (-1)^{xy} \{ <x|z>Py\  + \ <x|Pz>y
\} \eqno(2.8)$$
gives a quasi-classical $\delta$ Lie-super triple system.  Moreover, we 
have 
$$[Px,Py,Pz] = cP [x,y,z] \quad .$$
If $P = {1 \over 2} Id$, then this case reduces to the example 2.1.

\medskip

\noindent \underbar{\bf Example 2.3}

Let $L$ be a quasi-classical Lie-super algebra $(\delta =1)$.  If we 
introduce a triple product $[x,y,z]$ in $L$ by
$$[x,y,z] = [[x,y],z] \quad, $$
then $L$ becomes a quasi-classical Lie-super triple system with
 $\delta =1$. We may note that we then have
$$<[x,y,u]|v>\ =\ <[x,y]|[u,v]>$$
from which we can verify the validity of Eqs. (2.3).

\medskip

\noindent \underbar{\bf Remark 2.2}

We can calculate $<x|y>_1$ of the Proposition 2.2 for our various examples. 
First, the case of example 2.3 gives
$$<x|y>_1\  = \  {\rm Tr} (adx\  ady)$$
i.e., the Killing form of the Lie-super algebra $L$. On the other side, we 
calculate
$$<x|y>_1\  = (N_0 -1)<x|y>$$
and
$$<x|y>_1\  = ({\rm Tr}\ P) <x|y> \ +\  (N_0 -2) <x|Py>$$
for examples (2.1) and (2.2), respectively.  Here, we have set 
$$N_0 =\ {\rm Tr}\ 1 = \dim V_0 - \dim V_1 \quad .$$
However, we find Tr $L(x,y) =0$ for all cases in accordance with the Remark 
2.1. $\quad \bull$

Because of an intimate relationship between Lie-super algebra and Lie-super 
triple system for $\delta =1$, we will hereafter restrict ourselves to 
consideration only of the case $\delta =1$, unless it is stated otherwise.

\medskip

\noindent \underbar{\bf Remark 2.3}

Some connection exists between example 2.2 given above and example 1.1 or
 1.2 of the previous section. Let $L$ be the
 quasi-classical Lie or Lie-super
 algebra of either 1.1 or 1.2. Let $P\ \epsilon\ End\ L$ be defined by
$$Pf = e\quad , \quad Pe = Px_j = Py_j = 0 \quad (j=1,2,\dots,n)$$
which satisfies $P^2 = 0$ and $<Px|y>\ =\ <x|Py>.$

We can readily verify that $[[x,y],z]$ coincides with the expression $[x,y,z
]$ given by Eq. (2.8) of the example 2.2 for the same $<x|y>$. $\quad 
 \bull$

As we stated in example 2.3, we can construct a quasi-classical Lie-super 
triple system from a quasi-classical Lie-super algebra. The converse is 
also true as we will see below. To see it, we first define $M$ to be a 
linear span of the left multiplication operator $L(x,y)$ defined by Eq. 
(2.5a), i.e.,
$$M = \bigg\{Y |Y = \sum_{j,k} c_{jk} L (x_j,y_k)\bigg\} \eqno(2.9)$$
for  constants $c_{jk}$. Then, $M$ is a Lie-super algebra because of 
the lemma (2.1). A straightforward generalization of the well known 
canonical construction method enables us to go further as follows. Consider 
$$L_0 = V \oplus M \eqno(2.10)$$
for a Lie-super triple system $V$. We introduce a commutator in $L_0$, by
$$\eqalignno{&[x,y] = L(x,y)\ \epsilon\ M \quad ,&(2.10a)\cr
&[L(x,y),z] = - (-1)^{(x+y)z} [z,L(x,y)] = [x,y,z]
 \ \epsilon\ V \quad . &(2.10b)\cr}$$
Then, $L_0$ can be readily verified to be a Lie-super algebra for grading 
of
$$\sigma (L(x,y)) = \{ \sigma (x) + \sigma (y) \}\  ({\rm mod}\ 2) \quad .
 \eqno(2.11)$$
In order to make both $M$ and $L_0$ be quasi-classical, we introduce 
bilinear form in $M$ and $L_0$ by
$$\eqalignno{&<L(x,y)|L(u,v)>\ = \ < [x,y,u]|v>\ = -(-1)^{(u+v)y} <x|[ u,v,y]> 
&(2.12)\cr
&<L(x,y)|z>\ = \ <z|L(x,y)>\ =0 &(2.13)\cr}$$
in addition to $<x|y>$.

The second relation in Eq. (2.12)
 is the result of Proposition 2.1.  Note that Eq. (2.12) 
 is consistent with $L(x,y) = -(-1)^{xy} L
(y,x)$ and $L(u,v) = -(-1)^{uv} L(v,u)$. However, we have to verify its 
well-definedness, i.e. we have to verify the validity of 
$$< L(x^\prime ,y^\prime )|L(u,v)>\ =\ <L(x,y)|L(u,v)>\quad ,$$
for all $u,\ v\ \epsilon\ V$ whenever we have $L(x^\prime, y^\prime) = L(x,
y)$. This is trivially
correct, since we will have $<[x^\prime ,y^\prime, u]|v>\ =\ <[x,y,u]|v>$, 
if we note that $L(x^\prime , y^\prime) = L(x,y)$ implies $[x^\prime ,
y^\prime ,u] = [x,y,u]$ for any $u\ \epsilon\ V$.

\medskip

\noindent \underbar{\bf Proposition 2.3}

The Lie-super algebras $M$ and $L_0$ constructed canonically from a 
quasi-classical Lie-super triple system $V$ are quasi-classical.

\medskip

\noindent \underbar{\bf Proof}

First we will show that $<L(x,y)|L(u,v)>$ defined by Eq. (2.12) is 
non-degenerate. Suppose that we have 
$$<\sum_{j,k} c_{jk} L(x_j ,y_k)|L(u,v)>\ =0$$
for all $u,\ v\ \epsilon\ V$. This implies the validity of
$$<\sum_{j,k} c_{jk}[x_j,y_k,u]|v>\ =0 \quad .$$
Because of non-degeneracy of $<.|.>$, this leads to 
$$\sum_{j,k} c_{jk} [x_j, y_k, u] =0$$
or equivalently $\sum_{j,k} c_{jk} L(x_j ,y_k) =0$, proving the 
non-degeneracy. Next, we note
$$<L(x,y)|L(u,v)>\ =\ <[x,y,u]|v>\ =0\quad ,$$
unless $\sigma (x) + \sigma (y) + \sigma (u) + \sigma (v) = 0$ (mod 2) so 
that we find $<L(x,y)|L(u,v)>\ =0$ unless we have $\sigma(L(u,v)) = \sigma
 (L(x,y))$. Similarly, we find the validity of 
$$<L(x,y)|L(u,v)>\ =(-1)^{(u+v)(x+y)} <L(u,v)|L(x,y)>\quad .$$
Finally the proof for the validity of
$$<[L(x,y),L(z,w)]|L(u,v)>\ =\ <L(x,y)|[L(z,w),L(u,v)]> \eqno(2.14)$$
 goes as follows. In order to avoid unnecessary complications due 
to the sign factors $(-1)^{xy}$ etc., we will prove it only for non-super 
case. We can always supply sign factors for the super case to prove the 
same. Then, Eq. (2.14) is equivalent to
$$<[L(x,y),L(z,w)]|L(u,v)>\ =\ -<[L(u,v),L(z,w)]|L(x,y)>
 \quad . \eqno(2.14^\prime)$$
The left side of Eq. (2.14$^\prime$) is computed to be
$$\eqalign{<[L(x,y)&,L(z,w)]|L(u,v)>\cr
&=\ <L([x,y,z],w) + L (z,[x,y,w])|L(u,v)>\cr
&=\  <L([x,y,z],w) - L ([x,y,w],z)| L(u,v)>\cr
&= -<[x,y,z]|[u,v,w]>\ +\ <[x,y,w]|[u,v,z]>\quad .\cr}$$
If we interchange $x \leftrightarrow u$, and $y \leftrightarrow v$ in 
this expression, we find the validity of Eq. (2.14$^\prime$). This 
completes the proof, and the fact  that $L_0$ is quasi-classical also 
 can be 
similarly proved. $\quad \bull$

\medskip

\noindent \underbar{\bf Remark 2.3}

The canonical construction of an analogue of $L_0$ does \underbar{not} 
work for the case of $\delta =-1$.

\medskip

\noindent \underbar{\bf Def. 2.1}

A non-zero sub-vector space $B$ of a
  $\delta$ Lie-super triple system $V$ is called an ideal 
of $V$, if we have
$$[B,V,V] \subseteq B \quad .$$

\medskip

\noindent \underbar{\bf Proposition 2.4}

If $B$ is a ideal of a quasi-classical Lie-super triple system $V$
 $(\delta =1)$, then $L
(B,V)$ and $B \oplus L(B,V)$ are ideals of quasi-classical Lie-super 
algebras $M$ and $L_0$, respectively.

\medskip

\noindent \underbar{\bf Proof}

It is straightforward.

\medskip

\noindent \underbar{\bf Proposition 2.5}

Suppose that every ideal $B$ of a quasi-classical
 $\delta$ Lie-super triple system 
$V$ satisfies the condition 
$$[B,B,V] \not= 0 \quad .$$
Then, $V$ is a direct sum of simple ideals $B_j$:
$$V = B_1 \oplus B_2 \oplus \dots \oplus B_t \quad .$$
Moreover, we have

\medskip
\item{(i)} $<B_j |B_k>\ =0 \quad {\rm if}\quad j \not= k$
\smallskip

\item{(ii)} $[B_j ,B_k,V] = 0 \quad {\rm if}\quad j\not= k \quad .$

\medskip

\noindent \underbar{\bf Proof}

Let $B$ be a maximal ideal of $V$ and set
$$B^\prime = \ <x|<x|B>\ = 0 ,\  x\ \epsilon\ V> \quad .$$
Then, $B^\prime$ is a ideal of $V$, satisfying
$${\rm (i)}\quad <B|B^\prime>\ =0 \qquad {\rm (ii)}\quad [B,B^\prime, V] = 0
\qquad {\rm (iii)}\quad B \cap B^\prime = 0 \quad .$$
The fact that $B^\prime$ is an ideal of $V$ follows immediately from the 
Proposition 2.1, since
$$<[B^\prime, V,V]|B>\ =\ <B^\prime |[V,B,V]>\ =0 \quad .$$
Moreover,
$$<[B,B^\prime, V]|V>\ =\ <B|[V,V,B^\prime]>\ =0$$
also because of Eqs. (2.3b) and (2.1c). The non-degeneracy of $<.|.>$ 
then requires $[B,B^\prime, V]=0$. Next, set $A= B \cap B^\prime$. Suppose 
that $A \not= 0$. Then, $A$ is clearly an ideal of $V$.  However, $[A,A,V]=
0$ which is a contradiction with the hypothesis. Since $B$ is assumed to be 
maximal, these imply
$$V = B \oplus B^\prime \quad . $$
Moreover, $B$ and $B^\prime$ satisfy the same conditions as $V$. Hence, 
repeating the same arguments for $B$ and $B^\prime$, we reach at the 
conclusion of the Proposition.

\medskip

\noindent \underbar{\bf Remark 2.4}

It is plausible that $L_0$ corresponding to a simple quasi-classical
 Lie-super triple system will also be simple. However, the question will be 
discussed elsewhere.  Note that $M$ may be semi-simple (rather than being 
simple) even when $L_0$ is simple. See ref. [6] for such an example.

\medskip

\noindent \underbar{\bf Remark 2.5}

The special case $\delta =1$ in example 2.1 has been studied in [6] in 
connection with the para-statistics. It has been shown there that both $M$ 
and $L_0$ lead to simple Lie-super algebras of the type osp$(n|m)$ [1].  
For other examples, see also ref. [6]. $\quad \bull$

In ending this section, we would like to make some comments on
 Freudendal-Kantor triple systems,
 [7], since they are intimately connected with Lie-triple
 systems. Let $V$ be a $Z_2$-graded vector space with triple product 
$xyz$. If it satisfies
$$\eqalign{uv(xyz) =\  &(uvx)yz + \epsilon (-1)^{(u+v)x+uv} x(vuy)z\cr
&+ (-1)^{(u+v)(x+y)} xy(uvz)\cr}\eqno(2.15)$$
for $\epsilon = \pm 1$, $V$ is called a generalized Freudenthal-Kantor 
triple system. Especially, any Lie-super triple system is a generalized 
Freudenthal-Kantor triple system for $\epsilon =-1$ with $xyz = [x,y,z]$. 
On the other side, if we have
$$xyz = \delta (-1)^{xy+yz+zx} zyx \eqno(2.16)$$
with $\epsilon =-\delta$ in addition, it defines a
 $\delta$ Jordan-super triple system.  
Returning to the general case, we introduce a linear multiplication 
operator $K(.,.):V \otimes V \rightarrow End\ V$ by
$$K(x,y)z = (-1)^{yz} xzy - \delta (-1)^{x(y+z)} yzx \eqno(2.17)$$
for $\delta = \pm 1$. When we have identity
$$K(xyz,w) + (-1)^{z(x+y)} K(z,xyw) + \delta (-1)^{y(z+w)} K(x, K(z,w)y) = 
0 \quad , \eqno(2.18)$$
$V$ is called a $(\epsilon , \delta)$ Freudenthal-Kantor triple system [8]. 

The special case $K(x,y) =0$ with $\epsilon =-\delta$ will reproduce the
$\delta$ Jordan-super triple system. We can construct Lie-super triple 
systems out of $(\epsilon , \delta)$ Freudenthal-Kantor systems.  Here, we 
will present the following proposition.

\medskip

\noindent {\bf \underbar{Def. 2.2}}

Let $V$ be a $\delta$ Jordan-super triple system with bilinear
 non-degenerate form $<x|y>$ satisfying

\medskip

\item{(i)} $<x|y>\ = 0 \quad {\rm unless} \quad \sigma (x) = \sigma (y)$

\smallskip

\item{(ii)} $<y|x>\ = \delta (-1)^{xy} <x|y>$

\smallskip

\item{(iii)} $<xyu|v>\ =\ <x|yuv>\quad .$

\medskip

\noindent Then, $V$ is called quasi-classical.

\medskip

\noindent {\bf \underbar{Proposition 2.6}}

Let $V$ be a quasi-classical $\delta$ Jordan triple system. We introduce 
the left multiplication operation
$$L \ :\ V \otimes V \rightarrow End\ V$$
by
$$L(x,y)z = xyz \eqno(2.19)$$
with inner product
$$<L(x,y)|L(u,v)>\  =\  <xyu|v>\ =
\ <x|yuv> \quad .\eqno(2.20)$$
The resulting  Lie-super algebra given by
$$[L(u,v),L(x,y)] = L(uvx,y) - 
\delta (-1)^{(u+v)x+uv} L(x, vuy) \eqno(2.21)$$
is quasi-classical.

\medskip

\noindent \underbar{\bf Proof}

We first prove the validity of
$$<xyu|v>\ = (-1)^{(x+y)(u+v)} <uvx|y> \eqno(2.22)$$
since
$$\eqalign{<xyu|v>\ &= \delta (-1)^{xy+(x+y)u} <uyx|v>\ = \delta 
(-1)^{xy+(x+y)u} <u|yxv>\cr
&= \delta (-1)^{xy+(x+y)u} \cdot \delta (-1)^{v(x+y) +xy} <u|vxy>\cr
&= (-1)^{(x+y)(u+v)} <uvx|y> \quad .\cr}$$
We will have then
$$<L(u,v)|L(x,y)>\ = (-1)^{(u+v)(x+y)} <L(x,y)|L(u,v)> \quad .$$
It is easy to see then that it defines a non-degenerate super-symmetric
bilinear form.  Finally, the validity of
$$<[L(u,v),L(z,w)]|L(x,y)>\ =\ < L(u,v)|[L(z,w),L(x,y)]>$$
can be similarly shown just as in the proof of Eq. (2.14$^\prime$), if we 
note Eq. (2.21) and (2.22) to calculate
$$<[L(x,y),L(z,w)]|L(u,v)>\ =\ <xyz|wuv>\ -\  (-1)^{w(x+y+z)+z(u+v)} <wxy|uvz
>\quad . \quad \bull$$ 

\medskip

\noindent \underbar{\bf Proposition 2.7}

Let $xyz$ be a quasi-classical $\delta$ Jordan-super triple product. Then, 
$$[x,y,z] = xyz - 
\delta (-1)^{xy} yxz$$
defines a
 quasi-classical  $\delta$ Lie-super triple system.

\medskip

\noindent \underbar{\bf Proof}

It is straightforward.

\medskip

\noindent \underbar{\bf Example 2.4}

Suppose that $<x|y>$ and  $P\ \epsilon\ End\ V$ satisfy conditions
  of the example 2.2. 
Then, the product
$$xyz =\ <x|y>Pz\  + \ <x|Py>z\  + \ <y|Pz>x\  + \ <y|z> Px$$
defines a
 quasi-classical
  $\delta$ Jordan triple product. Further $[x,y,z]$ constructed in 
Proposition 2.7 reproduces the example 2.2. $\quad \bull$

\medskip

\noindent \underbar{\bf Example 2.5}

Let $L$ be a nilpotent Lie-super algebra of length at most 4, i.e., $L_5 =0
$. Especially, the examples 1.3 and 1.4 of section 1 satisfy the condition. 
For any two constants $c_1$ and $c_2$, we introduce a triple product by
$$xyz = c_1 [x,[y,z]] + c_2 [[x,y],z]$$
which defines a $(\epsilon , \delta)$ Freudenthal-Kantor system trivially. 
This is because we have
$$uv(xyz) = u(xyz)v = (xyz) uv = 0$$
in view of $L_5 =0$. Moreover, if we choose $c_1 = c_2$, it gives a
 quasi-classical
 Jordan-super triple system for $\delta = -\epsilon =1$.

\medskip

\noindent \underbar{\bf Example 2.6}

Let $<x|y>$ satisfy
$$<x|y>\ = -\epsilon (-1)^{xy} <y|x>\quad .$$
Moreover suppose that $P\ \epsilon\ End\ V$ obeys the condition
$$<P x|y>\ =\ <x|Py> \quad .$$
When   we set
$$xyz =\ <y|Pz>x \quad ,$$
we can verify the fact that it defines a $(\epsilon, \delta)$
 Freudenthal-Kantor triple system.

\medskip

\noindent {\bf 3. \underbar{Application to Yang-Baxter Equation}}

Let $R(\theta)$ be an element of $End\ (V \otimes V)$ for a parameter 
$\theta$ which is called the spectral parameter. We introduce $R_{jk}
(\theta)\ \epsilon\ End\  (V \otimes V \otimes V)$ for $j<k,\ j,k =1,2,3$ to 
be exactly like the operation of $R(\theta)$
 operating  only in $j$th and $k$th copies 
of $V$ in $V \otimes V \otimes V$. If we have
$$R_{12} (\theta) R_{13}(\theta^\prime) R_{23}(\theta^{\prime \prime}) =
 R_{23} (\theta^{\prime \prime}) R_{13} (\theta^\prime) R_{12}(\theta)
 \eqno(3.1)$$
for parameter $\theta ,\ \theta^\prime,$ and $\theta^{\prime \prime}$ 
satisfying
$$\theta^\prime = \theta + \theta^{\prime \prime} \quad , \eqno(3.2)$$
then the relation is called Yang-Baxter equation (e.g. see [9]).  Although 
we can generalize our result to the case of super space, we will consider 
here only non-super case for simplicity. Suppose that $V$ possesses a
 non-degenerate
 bilinear symmetric inner product $<.|.>$ so that we have $<y|x>\ 
=\ <x|y>$. We can then introduce ([10] and [11]) two $\theta$-dependent 
triple products $[x,y,z]_\theta$ and $[x,y,z]^*_\theta$ satisfying

\medskip

\item{(i)} $<x|[y,u,v]_\theta>\ =\ <y|[x,v,u]^*_\theta >$\hfill
(3.3)
\smallskip

\item{(ii)} $R(\theta) (x \otimes y) = \sum^N_{j=1} [e^j, y, x]^*_\theta
 \otimes e_j = \sum^N_{j=1} e_j \otimes [e^j ,x,y]_\theta \quad .$ \hfill 
(3.4)

\medskip

Here, $e_j$ and $e^j\ (j=1,2,\dots,N)$ are basis and its dual basis of $V$, 
respectively.  Then, the Yang-Baxter equation (hereafter abbreviated as
 YBE) can be rewritten as a triple product relation
$$\eqalign{\sum^N_{j=1} [&v,[u,e_j,z]_{\theta^\prime}, 
[e^j, x,y]_\theta ]^*_{\theta^{\prime \prime}}\cr
&= \sum^N_{j=1} [u,[v,e_j,x]^*_{\theta^\prime} ,[e^j ,
z,y]^*_{\theta^{\prime \prime}} ]_\theta \quad .\cr}\eqno(3.5)$$
We are hereafter interested only in the case when we have
$$[x,y,z]^*_\theta = [x,y,z]_\theta \eqno(3.6a)$$
or equivalently
$$<y|[x,v,u]_\theta>\ =\ <x|[y,u,v]_\theta> \quad . \eqno(3.6b)$$
Note that Eq. (3.6b) has the same form as Eq. (2.3c). Under these 
assumptions, we will first show:

\medskip

\noindent \underbar{\bf Lemma 3.1}

A necessary and sufficient condition to have
$$[R_{ij}(\theta),R_{k \ell} (\theta^\prime)] = 0 \eqno(3.7)$$
for all $i,j,k,\ell =1,2,3$ is the validity of
$$[u,v,[x,y,z]_\theta ]_{\theta^\prime} = [x,y,[u,v,z]_{\theta^\prime} 
]_\theta \quad . \eqno(3.8)$$

\medskip

\noindent \underbar{\bf Remark 3.1}

The validity of Eq. (3.7) implies that the YBE (3.1) as well as 
classical Yang-Baxter equation [9]
$$[R_{12}(\theta), R_{13}(\theta^\prime)] + [R_{12} (\theta),
 R_{23}(\theta^{\prime \prime}) ] + [R_{13}(\theta^\prime ), R_{23} 
(\theta^{\prime \prime}) ] = 0$$
hold valid without assuming the constraint Eq. (3.2).

\medskip

\noindent \underbar{\bf Proof}

We calculate for example 
$$\eqalign{R_{13}(\theta^\prime) R_{12}(\theta) x \otimes y \otimes z &= 
R_{13} (\theta^\prime) \sum^N_{j=1} [e^j ,y,x]_\theta \otimes e_j \otimes 
z\cr
&= \sum^N_{j,k=1} [e^k ,z,[e^j,y,x]_\theta]_{\theta^\prime} \otimes e_j 
\otimes e_k  \quad ,\cr
R_{12}(\theta) R_{13}(\theta^\prime) x \otimes y \otimes z &=
R_{12} (\theta) \sum^N_{k=1} [e^k ,z,x]_{\theta^\prime}
 \otimes y \otimes
e_k\cr
&= \sum^N_{j,k=1} [e^j ,y,[e^k,z,x]_{\theta^\prime}]_{\theta} \otimes e_j
\otimes e_k  \quad ,\cr}$$
from Eqs. (3.4) and (3.6a). Comparing both, we find 
$R_{12} (\theta) R_{13} (\theta^\prime) = R_{13} (\theta^\prime) R_{12}
(\theta)$ if we have Eq. (3.8). Similarly, we can prove the rest of 
relations. $\quad \bull$

\medskip

\noindent \underbar{\bf Lemma 3.2}

Let $L$ be a Lie algebra satisfying
$$[L,[[L,L],[L,L]]] = 0$$
as in the example 1.1. Then, the triple product defined by
$$[x,y,z] = [[x,y],z]$$
satisfies
$$[u,v,[x,y,z]] = [x,y,[u,v,z]]\eqno(3.9a)$$
or
$$[L(u,v),L(x,y)] = 0 \quad . \eqno(3.9b)$$

\medskip

\noindent \underbar{\bf Proof}

By a straightforward computation, we calculate
$$\eqalign{[u,v,[x,y,z]] &- [x,y,[u,v,z]]\cr
&= [[u,v],[[x,y],z]] - [[x,y],[[u,v],z]] \cr
&= [z,[[x,y],[u,v]]] = 0 \quad . \quad \bull\cr}$$

\medskip

\noindent \underbar{\bf Proposition 3.1}

Let $V$ be a quasi-classical Lie-triple systems satisfying
$$[u,v, [x,y,z]] = [x,y,[u,v,z]] \quad .$$
Then, $\theta$-dependent triple product
$$[x,y,z]_\theta = f(\theta) [x,y,z] + g (\theta )<x|y>z$$
for arbitrary functions $f(\theta)$ and $g(\theta)$ of $\theta$ gives a 
solution of Eq. (3.7), and hence of YBE.

\medskip

\noindent \underbar{\bf Proof}

The condition Eq. (3.6b) follows readily from Proposition (2.1), while we 
can easily verify the validity of Eq. (3.8). $\quad \bull$

\medskip

\noindent \underbar{\bf Proposition 3.2}

Let $L$ be a nilpotent quasi-classical Lie algebra of length at most 
4, i.e., $L_5 =0$. Then,
$$[x,y,z]_\theta = f_1 (\theta) [[x,y],z] + f_2 (\theta) [x,[y,z]] + g(
\theta) <x|y>z$$
for arbitrary functions $f_1(\theta),\ f_2(\theta)$ and $g(\theta)$ of 
$\theta$ is a solution of YBE.

\medskip

\vfill\eject

\noindent \underbar{\bf Proof}

If we set
$$<x,y,z>_\theta\  = f_1 (\theta) [[x,y],z] + f_2 (\theta) [x,[y,z]]\quad ,$$
it satisfies
$$<u,v,<x,y,z>_\theta >_{\theta^\prime}\  =\ <u, <x,y,z>_\theta ,
v>_{\theta^\prime}\  =\ <<x,yz>_\theta, u, v>_{\theta^\prime}\  = 0$$
as well as
$$<y|<x,v,u>_\theta >\ =\ <x|<y,u,v>_\theta > \quad .$$
Then, it is easy to check the validity of the required conditions
 Eqs. (3.6b) and (3.8). $\quad \bull$

\medskip

\noindent \underbar{\bf Remark 3.2}

Examples 1.3 of section 1 satisfies $L_4 =0$ and hence $L_5 =0$ of the 
condition. Another example satisfying Eq. (3.7) can be obtained as follows, 
although it does not correspond to a Lie triple system.  Let $J_\mu\ 
\epsilon\ End\ V$ for $\mu = 1,2,\dots ,m$ satisfy
$$[J_\mu , J_\nu] = 0 \quad .$$
Then,
$$R(\theta) = \sum^m_{\mu, \nu =1} f_{\mu \nu} (\theta) J_\mu \otimes
 J_\nu$$
for arbitrary functions $f_{\mu \nu} (\theta)$ of $\theta$ clearly satisfy 
Eq. (3.7). Such an example has been used elsewhere [12] to construct a 
rather curious link invariant.

\medskip

\noindent \underbar{\bf Remark 3.3}

We can find a solution of the YBE (3.5) but not necessarily of Eq. (3.7) as 
follows. Let $L$ be a nilpotent quasi-classical Lie algebra of length at 
most 6, i.e. $L_7 =0$. Then,
$$[x,y,z]_\theta = f_1 (\theta ) [x,[y,z]] + f_2 (\theta ) [[x,y],z] 
$$
is a solution of the YBE, which may not necessarily satisfy now Eq. (3.7).  
We can verify indeed that both sides of Eq. (3.5) vanish identically in 
view of $L_7 =0$. $\quad \bull$

\medskip

\noindent \underbar{\bf Acknowledgement}

The work of one of the authors (S.O.) is supported in part by U.S. 
Department of Energy Grant No. DE-FG02-91ER40685.

\vfill\eject

\noindent \underbar{\bf References}

\item{1.} M. Scheunert, \underbar{The Theory of Lie Superalgebras},
 (Springer-Verlag, Berlin 1979).

\item{2.} S. Okubo, A generalization of Hurwitz theorem and flexible 
 Lie-admissible algebra, Hadronic J. {\bf 3}, (1979) 1-52.

\item{3.} H.C. Myung, \underbar{Malcev-admissible Algebras}, (Birkh\"auser, 
Boston/Basel/Stuttgard, 1986).

\item{4.} R.D. Schafer, \underbar{An Introduction to Non-associative 
Algebras}, (Academic Press, New York/London, 1966).

\item{5.} N. Kamiya, A construction of anti-Lie triple systems from a class 
of triple systems, Mem. Fac. Sci. Shimane Univ. {\bf 22}, (1988) 51-62.

\item{6.} S. Okubo, Para-statistics as Lie-super triple systems, J. Math. 
Phys. {\bf 35}, (1994) 2785-2803.

\item{7.} N. Kamiya, A structure theory of Freudenthal-Kantor triple 
system, J. Alg. {\bf 110}, (1987) 108-123.

\item{8.} N. Kamiya, A structure theory of Freudenthal-Kantor triple system 
II, Comm. Math. Univ. Sancti Pauli {\bf 38}, (1989) 41-60.

\item{9.} M. Jimbo (ed.), \underbar{Yang-Baxter Equation in Integrable 
Systems}, (World Scientific, Singapore, 1989).

\item{10.} S. Okubo, Triple products and Yang-Baxter equation II. 
orthogonal and symplectic ternary systems, J. Math. Phys. {\bf 34}, (1993) 
3292-3315.

\item{11.} S. Okubo, Super-triple systems, normal and classical Yang-Baxter 
equation, in \underbar{Mathematics and
 its Applications} vol. 303, ed. by S. Gonzalez, 
(Kluwer Acad. Press, Dordrecht/Boston/London, 1994) pp. 300-308.

\item{12.} S. Okubo, New link invariants and Yang-Baxter equation, 
University of Rochester Report UR-1364 (1994).

\end